\documentclass[11pt]{article}
\usepackage{amsgen, amsmath, amsfonts, amsthm, amssymb, enumerate,amscd}

\topmargin -.15in \textheight 9in \numberwithin{equation}{section}

\newtheorem{defn}{Definition}[section]
\newtheorem{prop}[defn]{Proposition}
\newtheorem{lem}[defn]{Lemma}
\newtheorem{thm}[defn]{Theorem}
\newtheorem{cor}[defn]{Corollary}

\newcommand {\cf}{{\it cfr.}}
\newcommand {\ZZ}{{\mathbb Z}}

\newcommand {\K}{{\mathcal K}}

\newcommand {\C}{{\mathbb C}}

\newcommand {\F}{{\mathbb F}}

\newcommand {\Q}{{\mathbb Q}}
\newcommand {\R}{{\mathbb R}}
\newcommand {\OO}{{\mathcal O}}

\newcommand {\M}{{\mathcal M}}
\newcommand {\m}{{\mathfrak m}}

\newcommand {\PP}{{\mathfrak p}}

\newcommand {\CP}{{\mathbb P}}
\newcommand {\T}{{\cal T}}
\newcommand {\KI}{{K_{\infty}}}
\newcommand {\I}{{\mathcal I}}
\newcommand {\CI}{{{\mathbb C}_{\infty}}}

\title{$K_1$ of products of Drinfeld modular curves and
special values of $L$-functions}
\author{Caterina Consani and Ramesh  Sreekantan}
\date{}
\begin{document}
\sf \maketitle

\begin{abstract}

Let $X_0(I)$ be the Drinfeld's modular curve with level $I$ structure,
where $I$ is  a monic square-free ideal in  $\F_{q}[T]$. In this paper
we show  the existence of an  element in the  motivic cohomology group
$H^3_{\M}(X_0(I) \times X_0(I),\Q(2))$ whose regulator is related to a
special  value of  a Ranking-Selberg  convolution  $L$-function.  This
result is  the function field analogue  of a theorem  of Beilinson for
the self product of a modular curve.

\vspace{.1in}

{\bf Math. Subj. Class. Num.}: 11F52, 11G40.

\end{abstract}

\section{Introduction}

If $K$ is a real quadratic field with character $\chi$ of
conductor $N$, a classical formula in number-theory relates the
derivative at $s = 0$ of the $L$-function $L(s,\chi)$ to the
regulator of a unit in $K$:
\begin{equation}\label{cn1}
L'(0,\chi)=\log|U_{\chi}|.
\end{equation}
Here $U_{\chi}$ is the circular unit of $K$
\[
\prod_{k\text{mod}N\atop{(k,N) =1}}
(1-\zeta^{k})^{\frac{-\chi(k)}{2}}
\]
and  $\zeta=e^{2\pi i/N}$  is  a n-th  root  of unit.   There are  two
possible ways  of interpreting  \eqref{cn1}.  One is  in terms  of the
Beilinson's   conjectures  \cite{be}   relating   special  values   of
$L$-functions  of  algebraic varieties  over  number  fields to  their
$K$-groups. The  other interpretation is  in terms of a  conjecture of
Stark (\cf~\cite{stark}) on the  existence of certain units associated
to special values of  $L$-functions of number fields. Beilinson proved
his conjectures  in a few  cases - in  particular for elements  in the
$K_2$ group of  a modular curve and in the $K_1$  group of the product
of  two modular curves.   A closer  inspection of  Beilinson's theorem
shows that one  has an analogue of \eqref{cn1} for  the product of two
elliptic curves: \cf~\cite{ba-sr}.

In this  paper we are  concerned with a generalization  of
\eqref{cn1} for the  product of two Drinfeld  modular curves.  Our
main result is the following:

\begin{thm}
Let  $I=I_1 I_2$  be a  monic,  square-free ideal  in $\F_{q}[T]$  and
$X_0(I)$ be the Drinfeld modular  curve with level $I$ structure.  Let
$f$   and  $g$   be   two  normalized   new-forms   of  JLD-type   for
$\Gamma_0(I_1)$ and  $\Gamma_0(I_2)$ respectively. Then,  there exists
an element $\Xi_0(I)$ in the motivic cohomology group $H^3_{\M}(X_0(I)
\times X_0(I),\Q(2))$ such that
\begin{equation}\label{eq}
<reg(\Xi_0(I),f\bar{g})>~=~-\kappa
\frac{q-1}{\log_e(q)}~L'_{f,g}(0).
\end{equation}
Here, $\kappa \in \Q$ is  an explicitly determined number and
$L_{f,g}(s)$ is the Ranking-Selberg $L$-function convolution of
$f$ and $g$.
\end{thm}

In \cite{ko}, Kondo stated and  proved a similar statement for a
particular element in the $K_2$ group of  a Drinfeld modular
curve. In this paper we consider some of    his   ideas, together
with few results of    Gekeler
(\cf~\cite{ge1},\cite{ge2},\cite{ge3}, \cite{ge4}) and   Papikian
(\cf~\cite{pa}) to prove our result.

In the function field case, analogues of Beilinsons conjectures
are yet to be formulated. Our theorem, along with the results of
Kondo and Pal (\cf~\cite{pal}), provides some evidence for the
formulation of a general conjecture.

In the first part of the paper we introduce notation and
 describe the Bruhat-Tits tree of a Drinfeld modular
curve. It turns out that analogues of certain classical
number-theoretic quantities may be described in terms of this
combinatorial object.

In section~4 we recall the definition of an automorphic form of
JLD-type: this is the function field analogue of a classical
modular form of weight $2$.

Section~6 is dedicated to the computation of the  special value of
the L-function $L_{f,g}(s)$ associated to the convolution of two
automorphic forms $f$ and $g$ of JLD-type. Here, we use an
analogue of the classical Rankin-Selberg method that we describe
in section~5. The function $L_{f,g}(s)$ satisfies a functional
equation and we prove, using the theory developed by Gekeler, a
weak form of the Kronecker Limit formula.

Finally, section~7 contains a description of a relevant element in
the $K_1$ group of the self-product of a Drinfeld  modular curve.
We define its higher regulator and show that it is related to a
special value of $L_{f,g}(s)$. A consequence of our calculations
is an explicit proof of the finiteness of the cuspidal divisor
class group for modular curves of square-free level.

As a corollary of Theorem~1.1 we show the existence of an
interesting element in the motivic cohomology group of a product
of two elliptic curves of co-prime levels. Here, we use the
Taniyama-Shimura conjecture for function fields, which was proved
by Gekeler and Reversat in \cite{ge1}.\vspace{.05in}

As  our result  concerns  a  special value  of  the $L$-function,  the
constants that  appear in  the formula \eqref{eq}  might have  a finer
interpretation along the lines of the Bloch-Kato conjectures. It would
be interesting to investigate thoroughly this aspect.

\noindent{\bf  Acknowledgments.}  We have  greatly  profited from
the ideas contained in the preprint of S. Kondo \cite{ko} and in
the paper of M.  Papikian \cite{pa}.  We thank both  of them
warmly  for sending their papers to  us and for sharing some
insight  on this subject. The second author would also like to
thank J.Korman for some help.

The first author  is partially supported by the  NSERC grant
72016789; the second  author would like to  thank the University
of Toronto for its support.

\section{Notation}

In this paper we will use the following notation
\begin{itemize}
\item[{-}] $\F_q$ - the finite field with $q=p^{n}$ elements.

\item[{-}] $A=\F_{q}[T]$ - polynomial ring in one variable.

\item[{-}] $K=\F_q(T)$ - the quotient field of $A$

\item[{-}] $\pi=T^{-1}$ - a uniformizer at the infinite place
$\infty$

\item[{-}] $K_{\infty} = \F_q((\pi))$ - the completion of $K$ at
$\infty$

\item[{-}] $\CI$  - completed algebraic closure of $K_{\infty}$ -
the algebraic `complex plane'.

\item[{-}] $ord_{\infty} = - deg$  - the negative value of the
usual degree function.

\item[{-}] $\OO_{\infty} = \F_q[[\pi]]$ - $\infty$-adic integers.

\item[{-}] $|\cdot|$ - the $\infty$-adic absolute value on $\KI$
extended to $\CI$.

\item[{-}] $|\cdot|_{i}$ - the `imaginary part' of $|\cdot|$:
$|z|_{i}=inf_{x \in K_{\infty}}\{|z-x|\}$

\item[{-}] $G$ - group scheme $GL_2$.

\item[{-}] $B$ - Borel subgroup of $G$.

\item[{-}] $Z$ - Center of $G$.

\item[{-}] $\K = G(\OO_{\infty})$

\item[{-}] ${\cal I}$ - $\K(\pi)$ - Iwahori subgroup of $\K$.

\item[{-}] $\Gamma = G(A)$

\item[{-}] $\T$ - Bruhat-Tits tree of $PGL_2(K_{\infty})$.

\item[{-}] $V(\T)$ - the set of vertices $v$ of $\T$.

\item[{-}] $Y(\T)$ - the set of oriented edges of $\T$: if $e$ is
an edge, $o(e)$ and $t(e)$ denote the origin and terminus of the
edge.

\item[{-}] $\m$ - a divisor (ad\`ele) of $K$ with degree $deg(\m)$
(Note: This is different from $deg(m)=-ord_{\infty}(m)$ for $m \in
K$.)

\end{itemize}

\section{Preliminaries on the Bruhat-Tits tree}

The Bruhat-Tits tree $\T$ of $PGL_2(K_{\infty})$ is a
combinatorial object that is associated to a Drinfeld modular
curve. In this section we briefly recall the construction of the
tree $\T$. We refer to \cite {ge2} for more details.

\subsection {Vertices and Ends of $\T$}

Let $\T$ be the Bruhat-Tits  tree of $PGL_2(K_{\infty})$. This is
a tree whose set of vertices consist of classes of lattices $[L]$,
where $L$ is a $\OO_\infty$-lattice in $(K_{\infty})^2$ and $L$ is
said to be equivalent to $L'$ ($L \equiv L'$) if and only if there
exists an element $c \in K_{\infty}^*$ such that $L=cL'$. Two
vertices $[L]$ and $[L']$ are adjacent if they are represented by
lattices $L$ and $L'$ with $L \subset L'$ and
$\dim_{\F_{q}}(L'/L)=1$.  Every vertex $v$ has exactly
$q+1$-adjacent vertices and this set is isomorphic to
$\CP^1(\F_{q})$. More in general, the set of vertices of the tree
displayed at a distance $k$ from a fixed vertex $[L]$ is
isomorphic to $\CP^1(L/\pi^kL)$.

An end of $\T$ is an equivalence class of half-lines, with two
half  lines being equivalent  if and only if they  differ by  a
finite graph. Let $\partial \T$ be the set of the ends of $\T$.
There is a bijection (independent of $L$)
\[
\partial T \stackrel{\simeq}{\longrightarrow}\varprojlim_{k}\CP^1(L/\pi^kL) \simeq \CP^1(\OO_{\infty})
= \CP^1(K_{\infty}).
\]
The left-action of $G(K_{\infty})$ on $\T$ induces an action on
$\partial \T$ and this agrees with the action of $G(K_{\infty})$
on $\CP(K_{\infty})$ by linear fractional transformations.

\subsection {Orbit Spaces}

For $i\in\mathbb Z$, let $v_i \in V(\T)$ be the vertex
$[\pi^{-i}\OO_{\infty} \oplus \OO_{\infty}]$. Because the vertex
$v_0$ has stabilizer $\K\cdot Z(\KI)$ in $G(\KI)$, one obtains the
following identification
\[
G(K_{\infty})/\K\cdot Z(\KI) \stackrel{\simeq}{\rightarrow}
V(\T)\qquad g \mapsto g(v_0).
\]
Similarly, let $e_i$ be the edge $\overrightarrow{v_{i}v_{i+1}}$
({\it i.e}~$o(e_i) = v_i$, $t(e_i) = v_{i+1}$) then
\[
G(K_{\infty})/\I\cdot Z(\KI) \stackrel{\simeq}{\rightarrow}
Y(\T)\qquad g \mapsto g(e_0).
\]

These identifications  allows one to consider functions  on
vertices and on edges of $\T$ as equivariant functions on
matrices.

Let $w$  be the matrix $\begin{pmatrix} 0&1\\1&0 \end{pmatrix}$.
The following represents a standard choice of representatives for
$V(\T)$ and $Y(\T)$. Let
\[
S_{V}=\{\begin{pmatrix}\pi^k & u \\0  & 1 \end{pmatrix}|~k \in
\ZZ, u \in \KI, u \;\; \text{mod} \;\;\pi^k\OO_{\infty}\}.
\]

Then $S_V$ is  a system of representatives for $V(\T)$. Let
\[
S_U=\{w\begin{pmatrix}1 &0 \\c & 1\end{pmatrix} |~c \in \F_{q}\}
\cup \{1\},\quad S_Y=\{gh~|~g  \in  S_V,  h \in  S_U\}.
\]

Then, $S_Y$ is a system of representatives for $Y(\T)$. This
system  of representatives is  very important as any  function on
the vertices or edges of the tree $\T$ can be defined in terms of
it.

\subsection{Orientation}

The end  $\infty=(v_0,v_1,\ldots)$ defines  an orientation on $\T$
in the following manner. If $e$ is an edge, we say $e$ is positive
if it points  towards $\infty$  and negative if it points away
from infinity. This determines the decomposition
$Y(\T)=Y(\T)^{+}\cup Y(\T)^{-}$. Departing from any vertex $v$
there are precisely one positively oriented edge and $q$
negatively oriented edges with origin $v$. This determines a
bijection of $S_V$ with the set  of positively oriented  edges
$Y(\T)^+$. We will use the  notation $v(k,u)$ and $e(k,u)$ to
denote {\it resp.} the vertex and  the positively oriented edge
represented  by the matrix $\begin{pmatrix} \pi^k & u \\ 0  & 1
\end{pmatrix}$.  The edge  $e(k,u)$ has  origin $v(k,u)$  and terminus
$v(k-1,u)$.

\subsection{Realizations and norms}

The realization  $\T(\R)$  of  the  tree $\T$ is a topological
space that consists of  a real unit real interval for each
non-oriented edge glued together at the end  points according to
the incidence relations on $\T$.  If  $e$ is an oriented edge,
denote by $e(\R)$ the interval corresponding  to it  on the
realization. Let $\T(\ZZ)$ denote the points on $\T(\R)$
corresponding to the vertices of $\T$. The set of points
$\{t[L]+(1-t)[L']~|~t \in \Q\}$ lying on  edges $([L],[L'])$ will
be denoted by $\T(\Q)$.

A norm on  a  $\KI$-vector space  $W$  is  a function $\nu:W
\rightarrow \R$ satisfying the following properties
\begin{itemize}
\item[{-}] $\nu(v) \geq 0;~\nu(v)=0 \Leftrightarrow v=0$

\item[{-}] $\nu(xv)=|x|\nu(v),~\forall~x \in \KI$

\item[{-}] $\nu(v+w)\leq \text{max}\{\nu(v),\nu(w)\},~\forall~v,w
\in W$.
\end{itemize}

Two norms  are said to  be similar if  they differ by a non-zero
real constant.  The right action of $GL(W)$ on $W$ induces  an
action on the set of norms as
\[
\gamma(\nu)(v)=\nu(v\gamma).
\]
This action descents on similarity classes. The following theorem
will be crucial in what follows.

\begin{thm}[Goldman-Iwahori]
There is  a canonical  $G(\KI)$-equivariant bijection $b$ between
the set $\T(\R)$ and similarity classes of norms on $W=(\KI)^2$.
\end{thm}

The bijection $b$ is  defined as follows.  For  a  vertex  $v=[L]$
in $\T(\ZZ)=V(\T)$ we associate the norm $\nu_{L}$ defined by
\[
\nu_{L}(v) := \text{inf}\{|x|~|~x \in (\KI)^*, v \in xL\}.
\]
This norm makes  $L$ the unit  ball. If $P$  is a point of
$\T(\R)$ which lies on the  edge $([L],[L'])$ with $\pi L' \subset
L \subset L'$ and $P=(1-t)[L]+t[L']$, then $b(P)$ is the class of
the norm defined by
\[
\nu_P(v) := \text{max}\{\nu_{L}(v),q^t\nu_{L'}(v)\}.
\]

\section{The Drinfeld upper-half plane and the building map}

Let  $\Omega=\CP(\CI) -  \CP(\KI)=\CI -  \KI$ denote the Drinfeld
upper-half plane.  This space has  the structure of a rigid
analytic space over $\KI$. There is a canonical, surjective,
$G(\KI)$-equivariant map
\begin{equation}\label{bm}
\lambda:\Omega \longrightarrow \T(\R)
\end{equation}
called the building map. It is defined as follows.  To $z \in
\Omega$, we associate  the similarity  class of the  norm
$\nu_{z}$  on $(\KI)^2$ defined by
\[
\nu_{z}((u,v)) := |uz+v|.
\]
Since $|\;|$  takes values  in $q^{\Q}$, the  image of $\lambda$
in contained in $\T(\Q)$ and in fact $\lambda(\Omega)=\T(\Q)$. The
map $\lambda$ is used to describe an admissible covering
$\tilde\Omega$ of $\Omega$ (and of $\Gamma\backslash\Omega$). If
$v \in X(\T)$ and $e \in Y(\T)$: $\lambda^{-1}(v) \simeq \mathbb
P(\C_\infty)$ minus $(q+1)$-discs and $\lambda^{-1}(e) \simeq \{z
\in \C_\infty~|~|\pi|\le|z|\le 1\}$. Using $\lambda$ one obtains a
canonical identification of a certain intersection graph ({\it
i.e.} $(q+1)$-regular tree) of $\tilde\Omega$ with $\T$

\subsection{Drinfeld modular curves}

Let $\Gamma=G(A)$. For an ideal $I$ of $A$, we define
\[
\Gamma_0(I)=\{\begin{pmatrix} a & b \\
c & d \end{pmatrix}\in \Gamma~|~c \equiv 0 \;\; \text{mod}
\;\;I\}.
\]
Let $X_0(I)$  be the Drinfeld  modular curve of  level $I$. This
is a smooth, irreducible  algebraic curve  defined over a  finite
(abelian) extension  of $K$  such  that  there is  a  canonical
isomorphism  (as analytic spaces over $\CI$)
\[
X_0(I)(\CI) \simeq \Gamma_0(I)\backslash \Omega\cup \{cusps\}
\]
Entirely  analogous to  the corresponding  construction over  a
number field,  this curve  parametrizes  Drinfeld modules  of
rank two  with $I$-level  structure. The  corresponding quotient
of  the Bruhat-Tits tree will  be denoted by $\T_0(I)$.  Notice
that if $J  \subset I$ are two ideals of $A$, there is a
surjective map $X_0(J)\twoheadrightarrow X_0(I)$.

\subsection{Automorphic forms of JLD-type}

In the function field case, there  are two notions of modular
forms. One is the analogue of the classical  modular form on the
upper-half plane. The other, is the so-called automorphic form of
Jacquet-Langlands-Deligne (JLD) type. These are harmonic co-chains
on the edges of the tree $\T$. If $R$ is a commutative ring, a
$R$-valued harmonic cochain on $Y(\T)$ is a map $\phi:Y(\T)
\longrightarrow R$ satisfying the harmonic conditions:
\begin{itemize}
\item[{-}] $\phi(e) + \phi(\overline{e})=0$

\item[{-}] $\displaystyle{\sum_{t(e)=v}} \phi(e)=0.$
\end{itemize}

The second condition can also be stated as follows. Firstly,
notice that there is precisely one edge $e_0$  with $t(e_0)=v$ and
$sgn(e_0)=-1$ ({\it i.e.} $e_0\in Y(\T)^-$). The second condition
is then equivalent to
\[
\phi(e_0)=\sum_{t(e)=v\atop sgn(e)=1} \phi(e).
\]
If $\Gamma$ is an arithmetic subgroup of $G(A)$, we consider
$\Gamma$-invariant co-chains satisfying the further condition
\begin{itemize}
\item[{-}] $\phi(\gamma e)=\phi(e),\quad \forall \gamma \in
\Gamma.$
\end{itemize}

\noindent  The   group  of  $\Gamma$-invariant,   $R$-valued
harmonic co-chains  on edges is  denoted by
$H(Y(\T),R)^{\Gamma}$.The harmonic functions on  the edges are the
analogues of classical  cusp forms of weight  $2$.   In fact,  if
$\ell  \neq p$  is  a  prime number,  the $\Gamma$-invariant
harmonic co-chains  detect `half'  of  the \'etale cohomology
group $H^1(X(\Gamma),\Q_{\ell})$.

An automorphic  form of  level  $I$  is  an element $f \in
H^0(Y(\T_0(I)),\C) = H(Y(\T),\C)^{\Gamma_0(I)}$.   If  $f$ has
finite support, it is called a cusp form.

\subsection{Fourier expansions}

A  harmonic  function  on  the  set of positively oriented edges
$Y(\T)^+$ which is invariant under the group
\[
\Gamma_{\infty}=\{\begin{pmatrix} a & b \\0 & d \end{pmatrix} \in
G(A)\}
\]
has a Fourier expansion. This statement derives from the general
theory of Fourier analysis on ad\`ele groups: see \cite{ge2} for
details. The expansion has the following description. Let $\eta$
be the character $\eta:\KI \rightarrow {\mathbb C}^*$ defined as
\[
\eta(\sum_j a_j \pi^j) = \exp(\frac{2 \pi i Tr(a_1)}{p})
\]
where $Tr$ is the trace from  $\F_q$ to $\F_p$ (note: the role of
the two $\pi$'s in the formula is different!). Then, the Fourier
expansion of a $\Gamma_{\infty}$-invariant function $f$ on $Y(\T)$
is
\[
f(\begin{pmatrix} \pi^k & u \\0 & 1
\end{pmatrix}) = c_0(f,\pi^k) + \sum_{0 \neq m \in A \atop deg(m)
\leq k-2} c(f,div(m)\cdot \infty^{k-2}) \eta(mu)
\]
where the `constant' Fourier  coefficient $c_0(f,\pi^k)$ is the
function of $k \in \ZZ$ given by
\[
c_0(f,\pi^k) = \begin{cases} f(\begin{pmatrix} \pi^k & 0 \\0 & 1
\end{pmatrix}) &  if \; k\leq  0 \\ q^{1-k}  \displaystyle{\sum_{u \in
(\pi)/(\pi^k)} f(\begin{pmatrix}  \pi^k & u \\0 &  1
\end{pmatrix})} & if \; k \geq 1. \end{cases}
\]
For $\m$  a non-negative  divisor  on $K$  with $\m=div(m) \cdot
\infty^{deg(\m)}$, the `non-constant' Fourier coefficient is given
by
\[
c(f,\m)=q^{-1-deg(\m)} \sum_{u    \in (\pi)/(\pi^{2+deg(\m)})}
 f(\begin{pmatrix} \pi^{2+deg(\m)}& u \\0 & 1  \end{pmatrix})
\eta(mu).
\]

\subsection{Hecke operators and Hecke eigenforms}

Let $\PP=(p)$ be a prime ideal of $A$. The Hecke operator
$T_{\PP}$ is an operator on $H(Y(\T),\C)^{\Gamma_0(I)}$ defined by
\[
T_{\PP}(f)(e)=\begin{cases}  f(\begin{pmatrix} p & 0 \\ 0  & 1
\end{pmatrix} e )  + \displaystyle{\sum_{r~\text{mod}~\PP}} f(\begin{pmatrix}  1 & r
\\ 0 & p \end{pmatrix} e ) & if\;\;\PP \not{|}~I\\
\displaystyle{\sum_{r~\text{mod}~\PP}} f(\begin{pmatrix}  1 & r \\
0 & p
\end{pmatrix} e  ) & if \;\;\PP~|~I. \end{cases}
\]

A Hecke  eigenform is an automorphic form  of level $I$ which is
an eigenform for all  the Hecke operators  $T_{\PP}$.  Since the
automorphic forms are $\Gamma_{\infty}$-invariant,  they have
Fourier expansions.  The  Fourier coefficients of  cuspidal Hecke
eigenforms are known to have the following special properties
($\lambda_{\PP} \in \C$)
\begin{itemize}
\item[{-}] $c_0(f,\pi^k)=0,\quad \forall k \in \ZZ$

\item[{-}] $c(f,(1))=1$

\item[{-}] $c(f,\m)c(f,{\mathfrak n}) = c(f,\m {\mathfrak n})$,
whenever $\m$ and ${\mathfrak n}$ are relatively prime

\item[{-}] $c(f,\PP^{n-1})-\lambda_{\PP} c(f,\PP^n) +
|\PP|c(f,\PP^{n+1})=0$, when $\PP \not{|}~I \cdot \infty$

\item[{-}] $c(f,\PP^{n-1})-\lambda_{\PP} c(f,\PP^n)=0$ when
$\PP~|~ I$

\item[{-}] $c(f,\infty^{n-1})=q^{-n+1}$ if $n \geq 1.$
\end{itemize}


\subsection{Petersson inner product}

There is  an analogue  of the Petersson  inner product  for
automophic forms of JLD type. If $f$ and $g$ are automorphic forms
of level $I$, one of which is a cusp form, then we define
\[
\delta(f,g):=f(e)\bar{g}(e)d\mu(e)
\]
and the Petersson inner product
\[
<f,g>~:=   \int_{Y(\T_0(I))}   \delta(f,g)=   \int_{Y(\T_0(I))}
f(e) \bar{g}(e)d\mu(e)
\]
where $d\mu(e)$ is the Haar measure on the discrete set
$Y(\Gamma_0(I)\backslash\T)$ given by
$\frac{q-1}{2}|\text{Stab}_{\Gamma_0(I)}(e)|^{-1}$. Here,
$|\text{Stab}_{\Gamma_0(I)}(e)|$ denotes the cardinality of the
stabilizer of $e\in Y(\T)$.

If $f$  and $g$ are  normalized Hecke eigen-forms  and $f \neq g$,
then $<f,g>~= 0$.

\subsection{Logarithms and the logarithmic derivative}

If $f$ is  a function on a Drinfeld modular curve,  there is a
notion of a  logarithm defined  as follows. If  $v$ is a vertex on
$\T$  and  $\tau_{v}  \in  \Omega$ is an element of
$\lambda^{-1}(v)$ ({\cf~\eqref{bm}), define
\[
\log|f|(v):=\log_q|f(\tau_v)|.
\]
This quantity is  well defined as $|~|$  factors through the
building map. If $f$ is  a function on the vertices of the tree
$\T$, then the derivative of  $f$ is a function on the edges
defined by
\begin{equation}
\label{log}
\partial f(e): = f(t(e))-f(o(e)).
\end{equation}
The logarithmic derivative of a function $f$ is the composite of
these two, namely
\begin{equation}
\label{dlog}
\partial \log|f| (e) := \log|f|(t(e))-\log|f|(o(e)).
\end{equation}

\section{Eisenstein series and modular units}

The main goal of this section is the computation of a special
value of the  convolution $L$-function  of two  automorphic forms
$f$  and $g$ verifying  certain prescribed  conditions.   To  do
this we  study  some Eisenstein series  on the Bruhat-Tits tree
$\T$. We recall  that the classical Eisenstein-Kronecker-Lerch
series are  real   analytic functions on  a modular curve which
are related to  the logarithms of modular units  via the Kronecker
Limit formulas.   It turns  out that there  are function field
analogues  of these  series as well as an analogue of Kronecker's
First Limit formula. These results follow from the work of Gekeler
(\cf~\cite{ge2}) and are the crucial steps in the process of
relating the  regulators  of elements  in $K$-theory to special
values of $L$-functions.

\subsection{ Eisenstein series}

Let $I$ be an ideal in $A$ defined by a monic generator. The real
analytic Eisenstein series for $\Gamma_0(I)$  is defined as
\[
E_I(\tau,s)= \displaystyle{\sum_{\gamma \in
\Gamma_\infty\backslash\Gamma_0(I)}}
|\gamma(\tau)|_{i}^{-s},\qquad \tau \in \Omega,~s\in\C.
\]
This series converges for $Re(s) \gg 0$. The  function $|\;|_{i}$
factors  through  the building  map,  so  the Eisenstein series
can be thought of as a function defined on the vertices of the
Drinfeld modular curve. In terms of the matrix representatives
$S_{V}$, $E_I(\tau,s)$ is defined as follows.

Let $m,n \in A, (m,n)\neq (0,0)$ and let $v \in V(\T)$ be a vertex
represented  by  $v=\begin{pmatrix}  \pi^k  &   u\\0  &  1
\end{pmatrix}$. For $\omega=ord_{\infty}(mu+n)$ and $s\in\C$, define
\[
\phi_{m,n}^{s}(v)=\phi^s_{m,n} \begin{pmatrix} \pi^k  & u  \\ 0 &
1  \end{pmatrix} =
\begin{cases}  q^{(k-2deg(m))s}   &  if\;\;\  \omega   \geq  k-deg(m)\\
q^{(2\omega-k)s} & if \;\; \omega < k-deg(m). \end{cases}
\]
Then, using an explicit set of representatives for
$\Gamma_\infty\backslash\Gamma_0(I)$, we have (\cf~\cite{pa}
Section~4 for details)

\begin{equation}\label{E-I}
E_{I}(v,s)=q^{-ks}+\sum_{{m\in A\atop m\text{monic}}\atop m\equiv
0~\text{mod}I}\sum_{n \in A,\atop (m,n)=1} \phi^{s}_{m,n}(v).
\end{equation}

We define $E(v,s):=E_1(v,s)$.  In {\it loc.cit.},  it is shown
that $E_I(v,s)$ has an  analytic continuation to a meromorphic
function on the entire complex plane with a simple pole at $s=1$.
Lemma 3.4 of \cite{pa} relates the two series $E$ and $E_I$
through the formula
\[
\zeta_I(2s)E_I(v,s)=\frac{\zeta(2s)}{|I|^s}   \sum_{d~monic\atop
d|I} \frac{\mu(d)}{|d|^s} E((I/d)v,s)
\]
where $\zeta(s)=\frac{1}{1-q^{1-s}}$ is the zeta function,
$\mu(d)$ is the M\"obius function of $A$ and $(I/d)v$ denotes the
action of the matrix $\begin{pmatrix} I/d &0\\0 & 1 \end{pmatrix}$
on $v$.

Notice that a function  $G$ on the vertices of $\T$ can be
considered as a function on the  edges by defining
$G(e):=G(o(e))$. In particular, the function
\[
E_I(e,s):=E_I(o(e),s)
\]
agrees with the definition given in section~3 of \cite{pa}.

\subsubsection{Functional equation}

The  Eisenstein  series  $E(e,s)$  satisfies a  functional
equation analogous to the one verified in the classical case.
The ``Archimedean factor'' of the zeta function of $A$ is
$$L_{\infty}(s)=(1-|\infty|^{-s})^{-1}=\frac{1}{1-q^s}$$
Classically, the archimedean factor is obtained by multiplying the
$\Gamma$-function $\Gamma(s)$ by the factor $\pi^{-s}$. Here
instead, one multiplies $L_{\infty}(s)$ by $q^{s}$. In this way,
one gets $q^sL_{\infty}(s)=-\zeta(s+1)$
\begin{thm}
Let $\Lambda(e,s) := q^{s}L_{\infty}(s) E(e,s)$. Then
\[
\Lambda(e,s)=-\Lambda(e,1-s).
\]
Furthermore, $\Lambda(e,s)$ has a simple pole at $s=1$  with
residue $-(\log_{e} q)^{-1}$.
\end{thm}
\begin{proof} Cfr.~\cite{pa}~Theorem~3.3.
\end{proof}

\subsection{The Rankin-Selberg convolution}

An analogue of  the Rankin-Selberg formula for function field
automorphic forms may also be found in \cite{pa}. The Eisenstein
series $E_{I}(v,s)$ becomes an automorphic form by thinking of it
as a function on the edges as we have defined in the previous
section.

Let $f$ and $g$ be two automorphic forms of level $I$ on $\T$.
Consider the Dirichlet series
\[
L_{f,g}(s) = \zeta_I(2s) \sum_{\m \; \text{pos.div.}\atop
(\m,\infty)=1} \frac{c(f,\m) \bar c(g,\m)}{|\m|^{s-1}}
\]
where  $\zeta_I(s) = \prod_{\PP  \nmid I}  (1-|\PP|^{-s})^{-1}$. Since
$f$ and $g$ are normalized newforms, we have
\[
\zeta_I(2s)   \sum_{\m    \;   \text{pos.div.}}   \frac{c(f,\m)
\bar c(g,\m)}{|\m|^{s-1}}=-q^{s+1}L_{\infty}(s+1) L_{f,g}(s).
\]
In fact, following the decomposition $\m=\m_{fin}\infty^d$ ($d
\geq 0$), we obtain -from the last property of the Fourier
coefficients recalled in section~4.4-
$$c(f,\m)=c(f,\m_{fin}\cdot   \infty^d)=c(f,\m_{fin})q^{-d}.$$

Namely,  we can  pull out  the Euler factor at $\infty$.

\begin{prop}[Rankin's trick]
Let $f$ and $g$ be two cusp forms. Then
\[
\zeta_{I}(2s) <f~E_I(e,s),~g>~=~\zeta_I(2s) \int_{Y(\T_0(I))}
E_I(o(e),s)f(e)\bar g(e) d\mu(e) = -q^{2-s}L_{\infty}(s+1)L_{f,g}(s)
\]
\end{prop}

\begin{proof} \cf~\cite{pa}, section 4.
\end{proof}

This result determines an integral formula for $L_{f,g}(s)$. We
define
\begin{align*}
\Phi(s)&:=
\frac{\zeta_I(2s)q^{s}L_{\infty}(s)|I|^{s}}{\zeta(2s)}\int_{Y(\T_0(I))}
E_I(o(e),s) f(e) \bar g(e) d\mu(e) \\
&= \sum_{d~\text{monic}\atop d|I} \frac{\mu(d)}{|d|^s}
\int_{Y(\T_0(I))} \Lambda ((I/d)e,s)f(e)\bar g(e)d\mu(e) =
\frac{-q^{2}L_{\infty}(s)L_{\infty}(s+1)|I|^{s}}{\zeta(2s)}L_{f,g}(s).
\end{align*}
>From this we have the following observation
\begin{thm}[``Rankin's Theorem'']
$L_{f,g}(s)$ has a simple pole at $s=1$ with residue a multiple of
$<f,g>$. In particular, if $f$ and $g$ are normalized newforms and
$f \neq g$, then $L_{f,g}(s)$ is defined at $s=1$.
\end{thm}

We are interested in understanding the behavior of $L_{f,g}(s)$ at
$s=0$. As a first step we state the following theorem

\begin{thm} Let $f$ and $g$ be two cuspidal eigen-forms of square-free levels
$I_1,I_2$ respectively, with $I_1$ and $I_2$ co-prime ideals in
$A$. Let $I = I_1I_2$. Then, the function $\Phi(s)$ defined above
satisfies the functional equation
\begin{equation}
\Phi(s) = -\Phi(1-s). \label{fe}
\end{equation}

\end{thm}

\begin{proof}
The method of the proof is similar to that exposed in \cite{ogg}
(\cf~section 4). We show that by using the Atkin-Lehner operators,
we can simplify the integral. Let  $\PP|I$, say $\PP|I_1$.  We
have an  Atkin-Lehner  operator $W_{\PP}$ corresponding to $\PP$.
It can be represented by
\[
\beta=\begin{pmatrix} a\PP & -b \\I & \PP \end{pmatrix}, \;
\det(\beta) = \PP,\quad \beta \in \Gamma_0(I/\PP) \begin{pmatrix}
\PP & 0 \\0 & 1 \end{pmatrix}.
\]
where $a$ and $b$ are in $A$

Let $d~|~(I/\PP)$, then
\[
\begin{pmatrix} I / (\PP d) & 0 \\ 0 & 1
\end{pmatrix} \beta \begin{pmatrix} I/d & 0 \\ 0 & 1
\end{pmatrix}^{-1} \in \Gamma.
\]
Hence, since $\Lambda(e,s)$ is $\Gamma$-invariant one has
\begin{equation}
\Lambda((I/(d\PP))\beta e,s)=\Lambda((I/d)e,s). \label{beta}
\end{equation}
Since  $\beta$ normalizes $\Gamma_0(I_1)$ and $f$ is a new-form,
one obtains
\[
f|_{\beta}=f_{W_{\PP}}=c(f,\PP)f
\]
where $c(f,\PP)=\pm 1$. Further,  if $h=g|_{\beta}$ then
$h|_{\beta}=g|_{\beta^2}=g$. One then has
\begin{align}
\int_{Y(\T_0(I))} \Lambda((I/d\PP)e,s)\delta(f,g) &=
\int_{\beta^{-1}(Y(\T_0(I)))}
\Lambda((I/d)e,s)c(f,\PP)\delta(f,g|_{\beta}) \\
&= \int_{Y(\T_0(I))} \Lambda((I/d) e,s) c(f,\PP)\delta(f,h)\nonumber.
\end{align}
since $\beta^{-1}(Y(\T_0(I)))$ is a fundamental domain for $\beta^{-1}
\Gamma_0(I) \beta=\Gamma_0(I)$.

It follows that
\begin{equation}
\label{fe1} \Phi(s)= \sum_{d|(I/\PP)\atop
d~\text{monic}}\frac{\mu(d)}{|d|^s} \int_{Y(\T_0(I))}
\Lambda((I/d)e,s)( \delta(f,g)+c(f,\PP)\delta(f,h)).
\end{equation}
Now, we repeat this process with $h$ in the place of $g$. We first
observe  from the Fourier expansion,
\[
L_{f,h}(s)=c(f,\PP)|\PP|^{-s}L_{f,g}(s).
\]
Using this in the definition of $\Phi(s)$ we have
\begin{align}\label{fe2}
c(f,\PP)|\PP|^{-s}\Phi(s) &= \sum_{d|I/\PP\atop d~\text{monic}}
\int_{Y(\T_0(I))\atop d~\text{monic}}
(\Lambda((I/d)e,s)-|\PP|^{-s}\Lambda((I/(d\PP) e,s) ) \delta(f,h) \\
&= \sum_{d|I/\PP} \frac{\mu(d)}{|d|^{s}} \int_{Y(\T_0(I))}
\Lambda((I/d)e,s) \left( \delta(f,h)+\frac{c(f,\PP)}{|\PP|^{s}}\delta(f,g) \right).\nonumber
\end{align}
Let ${\mathfrak S}$ denotes the sum of the terms involving
$\delta(f,h)$. Comparing \eqref{fe1} and \eqref{fe2}, we obtain
$$ {\mathfrak S}(\frac{c(f,\PP)}{|\PP|^{s}}-1)=0.$$
Notice that the only way this equation  can hold for all $s\in\C$
is when ${\mathfrak S}=0$. Hence we have
\[
\Phi(s)=\sum_{d|(I/\PP)\atop d~\text{monic}}
\frac{\mu(d)}{|d|^{s}} \int_{Y(\T_0(I))} \Lambda((I/d)e,s)
\delta(f,g).
\]
Repeating this process for all primes $\PP$ dividing $I$, and
keeping in mind the assumption that a prime divides $I_1$ or
$I_2$ but not both, one gets
\begin{equation}
\Phi(s)=\int_{Y(\T_0(I))} \Lambda(Ie,s)\delta(f,g). \label{fe3}
\end{equation}
As $\Lambda(Ie,s)$ satisfies the functional equation
$\Lambda(Ie,s)= -\Lambda(Ie,1-s)$, we finally obtain
\[
\Phi(s)=-\Phi(1-s).
\]
\end{proof}

Hence, we conclude that if $f$ and $g$ are distinct eigen-forms:
$L_{f,g}(0)=0$,  as $L_{\infty}(s)$  has a  simple pole at  $s=0$.
This result is consistent with Beilinsons' conjectures  in the
number field case as one expects that the rank of  the motivic
cohomology is at least 1. In the number  field case, the  value of
the derivative of $L_{f,g}(s)$ at $s=0$ is related to the
regulator of an element of the motivic cohomology. To see whether
this relationship is verified in the function field case as well,
we need to get an explicit expression for $L'_{f,g}(0)$. This will
be accomplished in the next section.

\subsection{Kronecker's limit formula and the Delta function}

From the definition of $\Phi(s)$ given in the last section, we
have
\[
\Phi(0) = \frac{q^{2}}{\log_e(q)}L'_{f,g}(0).
\]
as  the  residue  of $L_{\infty}(s)$ at $s=0$ is
$\frac{-1}{\log_e(q)}$. For the computation of $\Phi(0)$ we need
to introduce the Drinfeld's $\Delta$ function.

\subsubsection{The Delta function}

Let $\tau$ be the coordinate function on $\Omega$ and let
$\Lambda_{\tau} = <1,\tau>$ be the  rank  two  free $A$-module in
$\CI$ generated by $1$  and $\tau$. Consider the following product
\[
e_{\Lambda_{\tau}}(z)= z \prod_{\lambda \in \Lambda_{\tau}
\backslash \{0\}} (1-\frac{z}{\lambda}) = z
\displaystyle{\prod_{a,b \in A\atop (a,b) \neq (0,0)}}
(1-\frac{z}{a \tau +b}).
\]
This product converges to give an entire, $\F_{q}$-linear,
surjective, $\Lambda_{\tau}$-periodic function
$e_{\Lambda_{\tau}}:\CI \rightarrow  \CI$.  This is the function
field  analogue of the classical $\wp$-function and it determines,
on the group scheme $\CI/\Lambda_{\tau}$, the structure of a
Drinfeld module.

Further, let $\Delta:\Omega \rightarrow \CI$ be the analytic
function defined by
\[
\Delta(\tau) = \mathop{\prod_{\alpha,\beta \in
T^{-1}A/A}}_{(\alpha,\beta) \neq  (0,0)} e_{\Lambda_{\tau}}(\alpha
\tau + \beta).
\]
This is the Drinfeld Delta  function. The Drinfeld's modular unit
$\Delta_I(\tau)$ is the function on the Drinfeld's upper half
space whose divisor is supported on the cusps:
\[
\Delta_{I}(\tau) = \prod_{d|I} \Delta(d\tau)^{\mu(\frac{I}{d})}.
\]

\subsubsection{The Kronecker Limit Formula}

In the classical case,  Kronecker's limit formula links the
Eisenstein series with  the logarithm  of the Delta  function. It
turns  out that there is an  analogue of this result in the
function field case. This follows from a theorem  of  Gekeler as
we are going to explain.

First of all we compute the constant term $a_0(v)$ in the Taylor
expansion of $E(v,s)$ around $s=1$. We have
\[
E(v,s)=\frac{a_{-1}}{s-1}+a_0(v)+a_1(v)(s-1)+\ldots
\]
where $a_{-1}$ is a constant independent of $v$. To compute
explicitly the coefficient-function $a_0(v)$ we differentiate
`with respect to $v$', namely we apply the $\partial$ operator
(\cf~section~4.6) and then we evaluate the result at $s=1$. This
gives
\[
\partial E(~,s)|_{s=1}=\partial a_0 (~).
\]
It follows that
\[
a_0(v)=\int \partial E(~,s)|_{s=1}(e) d\mu(e)+C
\]
where $C$ is a constant.

The function $\partial E(~,s)$ (on the edges $e\in Y(\T)$) is
related to the logarithmic derivative of the $\Delta$ function
through an improper Eisenstein series studied by Gekeler in
\cite{ge3}.

We first define Gekeler's series.  Let $\psi^s(e)= \text{sign}(e)
q^{-k(e)s}$. Consider the following Eisenstein series
\[
F(e,s)=\sum_{\gamma \in \Gamma_{\infty}\backslash \Gamma}
\psi^s(\gamma(e)).
\]
This series converges for $Re(s)\gg 0$. For $e=e(k,u) \in Y(\T)^+$
and
\[
\psi_{m,n}^{s}(e(k,u)) = \psi^s_{m,n}  \begin{pmatrix} \pi^k  & u
\\ 0  & 1  \end{pmatrix} =
\begin{cases}  -q^{(k-2deg(m)-1)s}   &  if\;\;\  \omega   \geq  k-deg(m)\\
q^{(2\omega-k)s} & if \;\; \omega < k-deg(m) \end{cases}
\]
where $\omega = \text{ord}_{\infty}(mu+n)$, we have
\[
F(e,s)=\psi^s(e)+ \sum_{{m \in A\atop m\text{monic}}\atop m\equiv
0~mod N}\sum_{n\in A\atop (m,n)=1} \psi^{s}_{m,n}(e).
\]

\begin{thm} The series $F(e,s)$  converges conditionally for $s=1$ and
\begin{equation}
\partial \log |\Delta|(e)=\frac{1-q}{q} F(e,1)
\label{gekeis}
\end{equation}
\end{thm}
\begin{proof} See \cite{ge3}, Corollary~2.8.\end{proof}

The following lemma determines a relation between the two series
$E(v,s)$ and $F(e,s)$.

\begin{lem}
For $E(v,s)$ and $F(e,s)$ defined as above,
\begin{equation}
\partial E(~,s)(e)=(q^s-1)F(e,s).
\label{partial}
\end{equation}
\end{lem}

\begin{proof}
\sf It follows from the definition of $\partial$ that
\[
\partial E(~,s)(e) = E(t(e),s) - E(o(e), s).
\]
Let $\omega = \text{ord}_{\infty}(mu+n)$ and
$e=e(k,u)=\stackrel{\longrightarrow}{v(k,u)v(k-1,u)}$. We will
consider four cases\vspace{.05in}

\noindent{\bf Case 0:} For $e=e(k,u)$
\begin{align}
\phi^s(t(e))-\phi^s(o(e)) &=q^{-(k-1)s}-q^{-ks}\\
&=(q^s-1)q^{-ks} \nonumber \\
&=(q^s-1)\psi^s(e) \nonumber
\end{align}

\noindent {\bf Case 1:} If $\omega > k-1-\text{deg}(m)$ then
\begin{align}
\phi^s_{m,n}(t(e))-\phi^s_{m,n}(o(e))
&= q^{(k-1-2deg(m))s}-q^{(k-2deg(m))s}  \\
&= (1-q^{s})q^{(k-1-2deg(m))s} \nonumber\\
&= (q^s- 1)\psi^s_{m,n}(e)\nonumber
\end{align}

\noindent {\bf Case 2:} If $\omega < k-1-\text{deg}(m)$ then
\begin{align}
\phi^s_{m,n}(t(e))-\phi^s_{m,n}(o(e)) &=
q^{(2\omega-(k-1))s}-q^{(2\omega-k)s}\\
&= (q^s-1)(q^{(2\omega-k)s}\nonumber \\
&= (q^s-1)\psi^s_{m,n}(e)\nonumber
\end{align}

\noindent {\bf Case 3:} If $\omega= k-1-\text{deg}(m)$ so
$2\omega-k=2k-2-2deg(m)$ then
\begin{align}
\phi^s_{m,n}(t(e))-\phi^s_{m,n}(o(e))
&= q^{(k-2deg(m)-1)s}-q^{(2\omega -k)s}\\
&= q^{(k-1-2deg(m))s}-q^{(2k-2-2deg(m))s}\nonumber \\
&= (q^s-1)(q^{(2\omega-k)s}\nonumber\\
&= (q^s-1)\psi^s_{m,n}(e).\nonumber
\end{align}
This shows that
\[
\partial E(~,s)(e)=(q^s-1)F(e,s).
\]
\end{proof}

Using this lemma we obtain the following result

\begin{thm}[``Kronecker's First Limit Formula'']
The function $\Lambda(v,s)$ has an expansion around $s=1$ of the
form
\begin{equation}
\Lambda(v,s)=\frac{b_{-1}}{s-1}+   \frac{q^2}{q-1}
\log|\Delta|(v)+ C+b_1(v)(s-1)+\ldots \label{klf1}
\end{equation}
where $b_{-1}$ and $C$ are constants independent of $v$.
\end{thm}

\begin{proof}

Recall   that   $\Lambda(v,s)=q^{s}L_{\infty}(s)E(v,s)$.   Therefore,
using \eqref{partial} we have
\[
\partial \Lambda(~,s)(e)=-q^{s}F(e,s).
\]
Evaluating this expression at $s=1$ and using \eqref{gekeis} we
obtain
\[
\partial\Lambda(~,1)(e)=\frac{q^2}{q-1}\partial
\log|\Delta|(e).
\]
Finally, by integrating we get the expansion of $\Lambda(v,s)$
around $s=1$
\[
\Lambda(v,s)=\frac{b_{-1}}{s-1}+   \frac{q^2}{q-1}
\log|\Delta|(v)+ C+b_1(v)(s-1)+\ldots
\]
\end{proof}

In \cite{ko}, Kondo proves the existence of a finer version of
this formula, {\it i.e.} an analogue  of Kronecker's second limit
formula. One can eventually use his result to obtain a finer
version of the next theorem.

\section{A special value of the $L$-function}

Using the functional equation for $\Lambda(v,s)$ and \eqref{klf1}
we obtain the following theorem

\begin{thm} Let $f$ and $g$ be two newforms with $f \neq
g$. Then
\begin{equation}
\Phi(0)= \frac{q^2}{q-1} <g \cdot \log|\Delta_I|,~f> =
\frac{-q^{2}}{\log_e q} L'_{f,g}(0). \label{spv}
\end{equation}

Here, $\log|\Delta_I|$ is thought of as a function on $Y(\T)$ as:
$\log|\Delta_I|(e)=\log|\Delta_I|(o(e))$.
\end{thm}

\begin{proof}
>From the definition of $\Phi(s)$ given in paragraph~5.2 we have
\[
\Phi(0)= \lim_{s \rightarrow 0} \sum_{d~\text{monic}\atop d|I}
\frac{\mu(d)}{|d|^{s}}\int_{Y(\T_0(I))}
\Lambda((I/d)e,s)\delta(f,g).
\]
Using the functional equation for $\Lambda(e,s)$ (\cf~Theorem~5.1)
and the Limit Formula \eqref{klf1}, we obtain

\begin{equation}
\Lambda(v,s)=\frac{b_{-1}}{s}-\frac{q^2}{q-1}\log|\Delta|(v)+C+
\text{h.o.t.(s)}\label{klf2}
\end{equation}

where $C$ is a constant. Therefore, using the fact that $<f,g>=0$
we have

\begin{align*}
\Phi(0) &= \lim_{s \rightarrow 0} \sum_{d~\text{monic}\atop d|I}
 \frac{\mu(d)}{|d|^{s}} \left( \int_{Y(\T_0(I))}  (\frac{b_{-1}}{s})
\delta(f,g)
-(\frac{q^2}{q-1}\log|\Delta(I/d)|(o(e))+C)\delta(f,g)+
\text{h.o.t}(s)\delta(f,g) \right) \\
& = \sum_{d~\text{monic}\atop d|I} \mu(d) b_{-1}+
\int_{Y(\T_0(I))} \sum_{d~\text{monic}\atop d|I} \mu(d)
\left(\frac{-q^2}{q-1} \right)\log|\Delta(I/d)|(o(e))\delta(f,g).
\end{align*}
The  term  involving $b_{-1}$  vanishes  as  $\sum_{d|I} \mu(d)=0$
since $I\neq 1$, so we have
\[
\Phi(0)=               -\frac{q^2}{q-1}
\int_{Y(\T_0(I))} \log|\Delta_I|(o(e))\delta(f,g).
\]
\end{proof}

\section{Elements in $K$-theory}

In  this section  we will  show  the the  special value of the
L-function computed in Theorem~\ref{spv} can be interpreted as the
regulator of an element in the $K$-theory group of the
self-product of Drinfeld modular curves. This statement should be
thought of as the function field analogue of the corresponding
theorem for products of modular curves proved by Beilinson  and it
supports some evidence for number field analogues of his
conjectures. We refer to \cite{ko} and \cite{pal} for a statement
in the case of $K_2$ of Drinfeld modular curves.

\subsection{The $K_1$ group of a surface}

If $S$  is an algebraic  surface, it is  well known that the
group of algebraic K-theory  $K_1(S)$ has a finite increasing
filtration - the Adams  filtration.  The second  graded  piece  is
usually denoted  by $H^3_{\M}(S,\Q(2))$:  this  is the  piece  of
the  filtration we  are interested in. When  $S$ is defined over a
number field, a conjecture of Beilinson  relates the co-volume of
the image of  this group under the  regulator  map into  a  real
vector space,  to  a  value of  the $L$-function of the middle
cohomology of the surface.

The group $H^3_{\M}(S,\Q(2))$ has the following description:
Elements in this group are represented by formal (finite) sums
\[
\sum_{i} (Z_i,u_i)
\]
where $Z_i$ are curves on $S$ and $u_i$ are functions on $ Z_i$
subject to the cocycle condition
\[
\sum_i div(u_i)=0.
\]
Relations  in   this  group   are  given  by   the  tame   symbols
of functions. More precisely, suppose that $Z$  is a curve on $S$
and $f$ and $g$ are two functions on $Z$. Then, the tame symbol is
defined by
\[
T_Z(f,g) = (-1)^{ord_Z(g)ord_Z(f)}
\frac{f^{ord_Z(g)}}{g^{ord_Z(f)}}.
\]
Therefore, elements like $\sum_Z (Z, T_Z(f,g))$ are zero in the
group $H^3_{\M}(S,\Q(2))$.

\subsection{Regulators}

Let $X_0(I)$  be the Drinfeld modular  curve of level $I$  an
ideal in $A$, and let $\Xi=\sum_i  (Z_i,u_i)$ be an element in
$H^3_{\M}(X_0(I) \times   X_0(I),\Q(2))$.   If  $f\overline{g}$ is
an  element   of $H^0(Y(\T_0(I)),\C) \times H^0(Y(\T_0(I)),\C)$,
the "$\infty$"-regulator is defined by
\[
<\text{reg}(\Xi),f\bar{g}> = <\text{reg}(\sum_i Z_i,u_i)),
f\bar{g}> = \sum_i \int_{Y(\T_{Z_i})} \log|u_i|(e) f(e)\bar
g(e)d\mu(e)
\]
where  $\T(Z)$ is the Bruhat-Tits tree corresponding to the curve
$Z$ on $X_0(I)\times X_0(I)$. Notice that the regulator is well
defined because of the Weil Reciprocity Law.

\subsection{A special element in $K$-theory}

In this paragraph we  will use the Drinfeld  modular  unit
\[
\Delta_I(\tau)=\prod_{d|I\;\text{monic}}
\Delta((I/d)\tau)^{\mu(d)}
\]
defined in section~5 along with  the diagonal of $X_0(I)$ to
construct a canonical element in the $K_1$ of the  self-product of
a Drinfeld curve. The trick we use here is  to `cancel  out'  the
zeroes and the poles of $\Delta_I$  using certain functions
supported on the vertical and  horizontal fibres of the variety.
The  existence of such functions is  a consequence of the analogue
of Manin-Drinfeld's theorem  which is proved  in the function
field case in \cite{ge3}.   Here, however,  we need  to work with
a quite explicit description  of these  functions. The goal is to
get an  effective version of the Manin-Drinfeld theorem.

\subsubsection{Cusps}

To compute  the divisor of $\Delta_I$ we  have to get a
handle of the cusps. The cusps of $X_0(I)$ are in bijection with
the set
\[
\{a/d~|~d~|~I~\text{monic},~a \in A/tA~\text{monic}\}
\]
where $t=(d,I/d)$. We will denote by $P^{a}_{d}$ the cusp
corresponding to $a/d$.

For simplicity, we will now  assume that $I$ is square-free.  Then
the cusps are  of the form $P_d=P^1_{d}$  where $d$ is a  monic
divisor of $I$. From now  on we will consider only monic divisors.
For a form or function F, we will use $F(f)$ to denote the form or
function $F(fz)$. For $a,b \in A$, let  $(a,b)$ denote the
greatest common divisor (gcd) of $a$ and $b$ and let  $[a,b]$
denote the least common divisor (lcm). We have the following
lemma.

\begin{lem}\label{thelem}
Assume that $I$ is a square-free, monic ideal of $A$ and that $I'$
and $d$ are monic divisors of $I$. Then
\begin{equation}
\text{ord}_{P_d} \Delta(I') = \rho_d~|I|~\frac{|(d,I')|}{|[d,I]|}
\label{ord}
\end{equation}
where $\rho_d=1$ if $d=1$ or $I$ and $\rho_d=(q-1)$ otherwise.
\end{lem}

\begin{proof}
\sf It follows from \cite{ge3} Lemma 3.8 that
\[
\text{ord}_{P_d}\Delta=\rho_d |I/d|,\qquad \text{ord}_{P_d}
\Delta(I)=\rho_d |d|.
\]
To compute the divisor of $\Delta(I')$ on $X_0(I)$ we need to
compute the ramification of $P_d$ over $P_{(d,I')}$, for some
divisor $d$ of $I$. Using {\em loc cit}, Lemma 3.8 we obtain
\[
\text{ram}^{P_d}_{P_{(d,I')}} = \frac{\rho_d
|I||(d,I)|}{\rho_{(d,I')}|d||I'|}.
\]
Therefore, it follows that
\begin{align*}
\text{ord}_{P_d} \Delta(I') &=  \text{ord}_{ P_{(d,I')}} \Delta(I')
\text{ram}^{P_d}_{P_{(d,I')}} = \frac{\rho_d
|I||(d,I)|}{\rho_{(d,I')}|d||I'|} \cdot
\rho_{(d,I)}|(d,I)| \\
&=  \rho_d |I| \frac{|(d,I')|}{|[d,I']|}.
\end{align*}

\end{proof}

   From  this  lemma we  deduce  that
\[  \text{div}(\Delta_I)=\sum_{d|I}
\mu(d)\text{div}(\Delta(I/d))= C \sum_{d|I} \mu(I/d) P_{I/d}
\]
where  $\mu(d)$ is  the  M\"obius function  on  $A$ and
$C=\sum_{d|I} \rho_d d $.\vspace{.1in}

We define a  simple modular unit to be a function  with divisor of
the form  $k(P)-k(Q)$  where  $P$  and  $Q$ are  cusps.  The
following theorem shows that there  is $\kappa \in \mathbb N$ such
that $\Delta_I^{\kappa}$ can be decomposed into a product of such
units.

\begin{thm}
\label{thethm} Let $I=\prod_{i=0}^{r} {f_i}$ be  the prime
factorization of $I$, where the $f_i's$ are monic in $A$. Let
$\kappa=\prod_{i=1}^{r} (|f_i|+1)$. Then, the function
$\Delta_I^{\kappa}$ decomposes as
\[
\Delta_I^{\kappa}=\prod_{d|\frac{I}{f_0}} F_{d}
\]
where the $F_d$'s  are simple units with divisor:
$\text{div}(F_d)=C \kappa(|f_0|-1)\mu(I/d)(P_d-P_{f_0d})$.

\end{thm}

\begin{proof}

\sf

Suppose that $I=f_0$. Then $\text{div}
(\Delta_I)=(|f_0|-1)(P_{f_0}-P_{1})$.  The general strategy of the
proof is to construct  a series of functions  with smaller and
smaller support until we will  obtain a function supported only on
two cusps.

Assume $f_1\neq f_0$ is a prime dividing $I$. Let consider the
following functions
\[
F_{f_1}(z)=\frac{\Delta(f_1z)^{|f_1|}}{\Delta(z)},\qquad
G_{f_1}(z)=\frac{\Delta(f_1 z)}{\Delta(z)^{|f_1|}}.
\]
In follows  from \eqref{ord} that
\[
\text{ord}_{P_d} F_{F_1}(f)=|f_1|\text{ord}_{P_d} \Delta(f) -
\text{ord}_{P_d} \Delta = \begin{cases} 0 \;& if\; (d,f_1)=1    \\
(|f_1|^2-1) ord_{P_d}\Delta\; &if\; f_1|d. \end{cases}
\]
A similar description holds for $G_{f_1}$. Let
\[
F=\prod_{d|(I/f_1)}F_{f_1}(dz)^{\mu(I/(f_1 d))};\quad
G=\prod_{d|(I/f_1)} G_{f_1}(dz)^{\mu(I/(f_1d))}.
\]
Then  $F$ and  $G$ are  functions supported  on complementary sets
of cusps and one has
\[
\Delta_I^{|f_1|+1}=FG.
\]

Continuing  in this manner  and replacing  $F_{f_1}$ for  $\Delta$
and $f_2$  for $f_1$,  we  can write  $F^{|f_2|+1}$  as a  product
of  two functions, one of which is supported on $P_d$ such that
$f_1f_2|d$ and the other  one on  the $P_d$ such  that $(f_2,d)=1,
f_1|d$.  A similar statement  holds  for  $G_{f_1}$.  As a result
we  can  write $\Delta_I^{(|f_1|+1)(|f_2|+1)}$ as product of four
functions supported on mutually exclusive  sets of cusps.
Repeating this  process for each $f_i|I$ with  $i \neq 0$,  we
obtain a  set of functions $F_d$  for each divisor  $d$   such
that  $div(F_d)$   is  supported  on   $P_d$  and $P_{f_0d}$. At
the  $i^{th}$ step one raises the  previous function to the
$|f_i|+1$-th power. This means that one has to raise $\Delta_I$ to
the $\kappa^{th}$ power.

\end{proof}

Notice that since the choice  of the prime $f_0$ was arbitrary,
there are  several different factorizations  of powers  of
$\Delta_I$. As a consequence of the computations shown before, we
obtain the following corollary of independent interest.

\begin{cor}[Manin-Drinfeld theorem]

Let $I=\prod_{i=0}^{r} f_i$ be the prime factorization of $I$.
Then, the cuspidal divisor class group is finite of order dividing
$C\prod_{i=0}^{r}(|f_i|^2-1)$
\end{cor}

\begin{proof}

If $P_{d_1}$  and $P_{d_2}$ are  two cusps, then $P_{d_1}-P_{d_2}$
can be written as  a sum of terms  of the form $P_{d}-P_{f_i  d}$
where $f_i$ are primes  dividing $I$. From  Theorem~7.2 follows
that the terms $P_d-P_{f_i}d$ are annihilated  by $C(|f_i|-1)
\prod_{j=0, j \neq i }^{r}(|f_j|+1)$. Hence  $P_{d_1}-P_{d_2}$
will be annihilated by the  least common multiple  of all these
numbers.

\end{proof}

\subsubsection{An element in $H^3_{\M}(X_0(I) \times X_0(I),\Q(2))$}

Let $D_0(I)$  denote the diagonal on $X_0(I)\times X_0(I)$ and let
$I=\prod_{i=0}^{r} f_i$ be the prime (monic) factorization of the
monic ideal $I \subset A$. Consider the element
\begin{equation}\label{el}
\Xi_0(I) := (D_0(I),\Delta_I^{\kappa})-\left( \sum_{d|(I/f_0)}
(P_{d} \times X_0(I),P_d \times F_d) - (X_0(I) \times P_{f_0d},
F_{d} \times P_{f_0d}) \right)
\end{equation}
It follows from Theorem~\ref{thethm} that this element satisfies a
cocycle condition as the divisor of it is a sum of multiples of
terms of the form
\[  (P_d,P_d)-(P_{f_0d},P_{f_0d})  -  (P_d,P_d)  +  (P_{d},P_{f_0d})  +
(P_{f_0d},P_{f_0d}) - (P_{d},P_{f_0d}).
\]
Hence $\Xi_0(I)$ determines an element of $H^3_{\M}(X_0(I) \times
X_0(I),\Q(2))$.

\section{The final result}

Our construction can be summarized by the following result

\begin{thm}
\label{mainthm}
If $I=I_1I_2$  is a monic,square-free  ideal of $A=\F_{q}[T]$  and $f$
and $g$  are two  Hecke newforms of  JLD type for  $\Gamma_0(I_1)$ and
$\Gamma_0(I_2)$ respectively and $f \neq  g$, then there is an element
$\Xi_0(I) \in H^3_{\M}(X_0(I) \times X_0(I),\Q(2))$ such that
\[
<reg(\Xi_0(I)),f\bar{g})> = -\kappa \frac{q-1}{\log_e(q)}
L'_{f,g}(0).
\]
Here, $L_{f,g}(s)$ denotes the Rankin-Selberg convolution of $f$
and $g$ and $\kappa = \prod_{i=1}^r(|f_i|+1)$ for $I =
\prod_{i=0}^rf_i$.
\end{thm}

\begin{proof}
Let  $\Xi_0(I)$  be  the  element defined in \eqref{el}. We
compute its regulator against the forms $f\bar{g}$
\begin{equation}
<\text{reg}(\Xi_0(I)),f\bar{g}> = \int_{Y(D_0(I))}\kappa
\log|\Delta_I(o(e))|f(e)\bar g(e) d\mu(e)= \kappa
<\log|\Delta_I|,f\bar{g}>. \label{reg}
\end{equation}

Note that the  regulator of the  vertical and horizontal
components computed against $f\bar{g}$  are $0$.  Combining
\eqref{reg} with Theorem~6.1, we obtain
\[
<reg(\Xi_0(I)),f\bar{g})> = -\kappa \frac{q-1}{\log_e(q)}
L'_{f,g}(0).
\]
\end{proof}

\subsection{Application to elliptic curves}

If $E$ is a non-isotrivial ({\it i.e.} $j_E \notin \F_q$)
semi-stable elliptic curve over $K$ with  split  multiplicative
reduction  at $\infty$ and of conductor $I_{E}=I\cdot\infty$,  it
is shown in \cite{ge1} that $E$  is modular. This means that the
Hasse-Weil Zeta  function  $L(E,s)$ is equal  to the  $L$-function
of an  automorphic  form $f$ of  JLD type with rational fourier
coefficients
\[
L(E,s)=L(f,s)=\sum_{\m \;\text{pos. div}}
\frac{c(f,\m)}{|\m|^{s-1}}.
\]
Furthermore, there exists a non-trivial morphism
\[
\pi:X_0(I) \longrightarrow E
\]
where $X_0(I)$ is the Drinfeld modular curve of level $I$.

Now suppose  $E$   and  $E'$  are  two  such   elliptic  curves
with corresponding  automorphic  forms $f$  and  $g$  of  levels
$I_1$ and $I_2$. Assume  that $(I_1,I_2)=1$ and that $I=I_1 I_2$
is a square-free ideal.  Then, the $L$-function of  $H^2(E \times
E')$ can be expressed in terms  of the $L$-function of  the
Rankin-Selberg convolution of $f$  and $g$. The K\"unneth theorem
gives the decomposition
\[
L(H^2(E    \times    E'),s)=    L(    H^2(E),s)L(H^1(E)\otimes
H^1(E'),s)L(H^2(E'),s)=\zeta(s-1)^2 L(H^1(E)\otimes H^1(E'),s)
\]
It is  not completely   clear to us what are the  terms
corresponding to the primes dividing $I$ in the expression for
$L(H^1(E)\otimes H^1(E'),s)$. On the other hand, the nice
functional equation for $\Phi(s)$ would suggest that the completed
$L$-function for $H^1(E)\otimes H^1(E)$ should be $\Phi(s-1)$. The
`Archimedean' term of $\Phi(s-1)$ is
$L_{\infty}(s-1)L_{\infty}(s)$ so we  divide out by it. Under this
assumption, we have
\begin{thm}
Let $E$ and $E'$ be elliptic curves over $K$ satisfying all the
above conditions. Then, there is an element $\Xi \in
H^3_{\M}(E\times E',\Q(2))$ such that
\[
L'(H^2(E \times E'),1)=\frac{\log_e(q)
deg(\Pi)}{(q-1)\kappa}<reg(\Xi),f \overline{g}>
\]
\end{thm}
\begin{proof}
\sf Let $\pi_f \times \pi_g: X_0(I) \times X_0(I) \rightarrow E
\times E'$ be  the product  of the  modular parametrizations. Let
$\Xi$  be the pushforward $(\pi_f  \times \pi_g)_*(\Xi_0(I))$ to
$H^3_{\M}(E \times E',\Q(2))$. Let $\Pi$  be the restriction of
$\pi_f  \times \pi_g$ to the diagonal  $D_0(I)$. This contributes
a factor  deg($\Pi$) to the equation. Then, the result follows
from Thereom~\ref{mainthm}.

\end{proof}

\vspace{.2in}

\noindent {\bf Caterina Consani}, Department of Mathematics,
University of Toronto, Canada.

\noindent email: kc\@@math.toronto.edu

\vskip .1in

\noindent {\bf Ramesh Sreekantan}, Department of Mathematics,
University of Toronto, Canada.

\noindent email: ramesh\@@math.toronto.edu


\begin{thebibliography}{[G-K-Z]}

\bibitem [Ba-Sr]{ba-sr}  Baba, Srinath and Sreekantan,  Ramesh {\em An
analogue of circular units  for products of elliptic curves},
preprint 2002.

\bibitem   [Be]{be} Be\u\i  linson,  A.A. {\em  Higher regulators  and
    values  of   $L$-functions.}    (Russian)  Current    problems  in
  mathematics, Vol.  24, 181--238, Itogi Nauki i Tekhniki, Akad.  Nauk
  SSSR,  Vsesoyuz. Inst.   Nauchn. i  Tekhn.  Inform., Moscow,   1984.


\bibitem [Ge2]{ge2} Gekeler,  Ernst-Ulrich Improper {\em Eisenstein series
on Bruhat-Tits trees.} Manuscripta Math. 86 (1995), no. 3,
367--391.

\bibitem [Ge-Re]{ge1}  Gekeler, E.-U.;  Reversat, M. {\em  Jacobians of
Drinfeld modular curves.} J. Reine Angew. Math. 476 (1996),
27--93.


\bibitem   [Ge3]{ge3}   Gekeler,   Ernst-Ulrich  {\em On   the   Drinfeld
discriminant function.} Compositio Math. 106 (1997), no. 2,
181--202.

\bibitem  [Ge4]{ge4} Gekeler,  Ernst-Ulrich {\em On  the  cuspidal divisor
class group of a Drinfeld modular curve.} Doc. Math. 2 (1997),
351--374


\bibitem [Ko]{ko} Kondo, Satoshi {\em Euler systems on Drinfeld Modular Curves and zeta values}  UTMS Preprint, 2002-6.


\bibitem [Mi]{mi}Mildenhall, Stephen J. M. {\em Cycles in a product of
    elliptic curves, and a group analogous to  the class group.}  Duke
  Math. J. 67 (1992), no. 2, 387--406.

\bibitem  [Og]{ogg} Ogg,  A. P. {\em On   a convolution  of $L$-series.}
  Invent. Math.7 1969 297--312.

\bibitem [Pal]{pal} Pal, Ambrus {\em personal communication}

\bibitem  [Pa]{pa} Papikian,  Mihran  {\em On  the  degree of  modular
parametrizations  over function  fields} J.  Number Theory  97
(2002), 317--349.

\bibitem   [St]{stark}   Stark,    H.   M.   {\em   $L$-functions   at
$s=1$. II. Artin $L$-functions  with rational characters.}
Advances in Math. 17 (1975), no. 1, 60--92.



\end{thebibliography}
\end{document}